\documentclass[a4paper]{article}

\usepackage{amsmath,amstext,amsbsy,amssymb,graphicx,subfig,mathdots,tikz,todonotes}

\usepackage{mathrsfs}
\usepackage{fullpage}
\usepackage[utf8]{inputenc}
\newcommand{\N}{\mathbb{N}}
\newcommand{\R}{\mathbb{R}}
\newcommand{\C}{\mathbb{C}}
\newcommand{\Rn}{\R^n}
\newcommand{\Cnn}{\C^{n\times n}}
\newcommand{\Cmn}{\C^{m\times n}}

\newcommand{\norm}[1]{\|#1\|}
\DeclareMathOperator{\sym}{sym_{_*}\!}
\DeclareMathOperator{\skw}{skew_{_*}\!}
\DeclareMathOperator{\tr}{tr}
\newcommand{\trik}{\tr_i^{k}}
\newcommand{\trok}[1][k]{\tr_1^{(#1)}}
\newcommand{\comp}{^{(k)}}
\newcommand{\inv}{^{-1}}
\DeclareMathOperator{\diag}{diag}
\newcommand{\phii}{\varphi}
\newcommand{\my}{\mu}
\newcommand{\A}{Z} 

\DeclareMathOperator{\GL}{GL}
\DeclareMathOperator{\U}{U}
\newcommand{\id} {I}
\renewcommand{\Re}{\mathop{\mathfrak{Re}} }
\renewcommand{\Im}{\mathop{\mathfrak{Im} }}
\newcommand{\abs}[1]{|#1|}
\DeclareMathOperator{\Log}{Log}
\newcommand{\matr}[1]{\begin{pmatrix}#1\end{pmatrix}}
\newcommand{\sort}{^\downarrow}

\newtheorem{Theorem}{Theorem}[section]
\newtheorem{ex}[Theorem]{Example}
\newtheorem{Proposition}[Theorem]{Proposition}

\begin{document}

\title{The minimization of matrix logarithms - on a fundamental property of the 
unitary polar factor}

\author{Johannes Lankeit%
\thanks{Johannes Lankeit, Fakult\"at f\"ur  Mathematik, Universit\"at Duisburg-Essen, Campus Essen, Thea-Leymann Str. 9, 45127 Essen, Germany, email: johannes.lankeit@uni-due.de}\; and 
Patrizio Neff%
\thanks{Patrizio Neff, Lehrstuhl f\"ur Nichtlineare Analysis und Modellierung, Fakult\"at f\"ur  Mathematik, Universit\"at Duisburg-Essen, Campus Essen, Thea-Leymann Str. 9, 45127 Essen, Germany, email: patrizio.neff@uni-due.de, Tel.: +49-201-183-4243}\; and 
Yuji Nakatsukasa%
\thanks{ Yuji Nakatsukasa, 
Department of Mathematical Informatics,  University of Tokyo, Tokyo 113-8656, Japan, 
email: nakatsukasa@mist.i.u-tokyo.ac.jp. Supported in part by EPSRC grant EP/I005293/1.}
}

 \maketitle


\begin{abstract}
 We show that the unitary factor $U_p$ in the polar decomposition of a nonsingular matrix $Z=U_pH$ is the minimizer for both 
\[ \norm{\Log(Q^*Z)}\quad\text{ and }\quad \norm{\sym(\Log(Q^*Z))}\]
 over $Q\in U(n)$ for any given invertible matrix $Z\in\Cnn$, for any unitarily invariant norm and any $n$. We prove that $U_p$ is the unique matrix with this property. 
As important tools we use a generalized Bernstein trace inequality and the theory of majorization.
\end{abstract}



{\bf keyword}
unitary polar factor, matrix logarithm, matrix exponential, Hermitian part, minimization, unitarily invariant norm, polar decomposition,
majorization, optimality



\section{Introduction}
Just as every nonzero complex number $z=re^{i\phii}$ admits a unique polar representation with $r\in\R_+, \phii\in (-\pi,\pi]$, every matrix $Z\in\Cnn$ can be decomposed into a product of the unitary polar factor $U_p\in U(n)$ (where $U(n)$ denotes the group of $n\times n$  unitary matrices) and a positive semidefinite matrix $H$ \cite[Lemma 2, p.124]{Autonne1902}, \cite[Ch.~8]{Higham2008},\cite[p.414]{Horn85}:
\begin{align*}
 \label{eq:polardecomposition}
 Z=U_p\, H.
\end{align*}
This decomposition is unique if $Z$ has full column rank. We note that the polar decomposition exists for rectangular matrices $Z\in\Cmn$, but in this paper we shall restrict ourselves to invertible $Z\in\Cnn$, in which case $U_p,H$ are unique and $H=\sqrt{Z^*Z}$ is positive definite, where the matrix square root is taken to be the principal one~\cite[Ch.~6]{Higham2008}.

The unitary polar factor $U_p$ plays an important role in geometrically exact descriptions of solid materials. In this case  $U_p^TF=H$ is called the right stretch tensor of the deformation gradient $F$ and serves as a basic measure of the elastic deformation \cite{NeffBirsan13,Neff_Chelminski_ifb07,Neff_Forest_jel05,Neff_micromorphic_rse_05,Neff01c}. For additional applications and computational issues of the polar decomposition see e.g. \cite[Ch.~12]{Golubbookori} and \cite{nakahig12, ByXu08,nakatsukasa:2700,nahi11}.

The unitary polar factor also has the property that in terms of any unitarily invariant matrix norm $\norm{\cdot}$, i.e. norms that satisfy $\norm{X}=\norm{UXV}$ for any unitary $U,V$, 
 it is the nearest unitary matrix \cite[Thm.~IX.7.2]{Bhatia97}, \cite{Fan55}, \cite[p.~197]{Higham2008} to $Z$, that is, 
\begin{equation}
 \label{irgendwas}
 \min_{Q\in U(n)} \norm{Z-Q}=\min_{Q\in U(n)} \norm{Q^*Z-I}=\norm{U_p^*Z-I}=\norm{\sqrt{Z^*Z}-\id}.
\end{equation}
The purpose of this paper is to
 show that the unitary polar factor enjoys this minimization property (made precise in \eqref{min_defi}) also with respect to $\norm{\Log \ldots}$, an expression that arises when considering geodesic distances on matrix Lie groups (see \cite{NNF}, \cite{NeffEidelOsterbrinkMartin_Riemannianapproach} and \cite{Neff_Osterbrink_Martin_hencky13} for further motivation):
\[
 \min_{Q\in U(n)} \norm{\Log Q^*Z}=\norm{\log U_p^*Z}=\norm{\log\sqrt{Z^*Z}},
\]
and with respect to the Hermitian part of the logarithm
\[
 \min_{Q\in U(n)} \norm{\sym\Log Q^*Z}=\norm{\sym\log U_p^*Z}=\norm{\log\sqrt{Z^*Z}}. 
\]
Here $\Log Z$ denotes any solution to $\exp X=Z$, while $\log Z$ denotes the principal matrix logarithm (we discuss more details in section~\ref{sec:matlog});
 $\sym X=\frac{1}{2}(X+X^*)$ is the Hermitian part of $X\in\Cnn$.

This minimization property is fundamental as it holds for arbitrary $n\in\N$, all unitarily invariant matrix norms, and in fact for the whole family
\begin{align}
\label{eq:minimization_family}
 \my\norm{\sym\Log(Q^*Z)}+\my_c\norm{\skw\Log(Q^*Z)}, \quad \my>0,\my_c\geq 0.
\end{align}
By contrast, the respective property does not hold true 
\cite{Neff_Biot07} for 
\begin{align}
 \label{eq:minimization_family_false}
 \my\norm{\sym(Q^*Z-\id)}+\my_c\norm{\skw(Q^*Z-\id)}, \quad 0<\my_c<\my,
\end{align}
wherefore the minimization \eqref{eq:minimization_family} seems even more fundamental than \eqref{eq:minimization_family_false}. 
Note that 
\eqref{eq:minimization_family_false} reduces to \eqref{irgendwas} by taking $\my=\my_c=1$.

This result, which is a generalization of the fact for scalars that for any complex logarithm and for all $z\in\C\setminus\{0\}$
\begin{align*}
 \min_{\vartheta\in(-\pi,\pi]}\abs{\Log_{\C}(e^{-i\vartheta} z)}^2&=\abs{\log\abs{z}}^2\, ,\quad
 \min_{\vartheta\in(-\pi,\pi]}\abs{\Re\Log_{\C}(e^{-i\vartheta} z)}^2=\abs{\log\abs{z}}^2\, ,
\end{align*}
has recently been proven for 
the spectral norm in any dimension $n$ and the Frobenius norm 
for $n\leq 3$
 in \cite{NNF}. 
By using majorization techniques (see also \cite{Bhatia2012726}) we now 
prove this property
\emph{in any dimension} $n$ 
and for \emph{any} unitarily invariant matrix norm.

In \cite{NNF} the conditions for applying the new sum of squared logarithms inequality \cite{Neff_log_inequality13} are obtained from the inequality 
\begin{align}
\label{eq:expleqexpsym}
\norm{\exp X}\leq \norm{\exp \sym X}, \qquad \text{\cite[IX.3.1]{Bhatia97}}, 
\end{align}
which can be derived from Cohen's generalization \cite{Cohen88} of Bernstein's trace inequality \cite{Bernstein88}, which is inequality \eqref{eq:expleqexpsym} for the Frobenius norm.
In this paper, we exploit the conditions obtained by Cohen \cite{Cohen88}, inequality \ref{eq:Cohen} below, directly, apply the logarithm first and then use majorization techniques.

In the next section we provide some basics about compound matrices and majorization upon which our proof is built.
We then discuss properties of the matrix logarithm, 
and in section \ref{sec:main} we prove the asserted minimization property. Finally, we prove the uniqueness of $U_p$ as the minimizer. 

\emph{Notation}. 
$\sigma_i(X)=\sqrt{\lambda_i(X^*X)}$ denotes the $i$-th largest singular value of $X$. 
The symbol $\id_k$ denotes the $k\times k$ identity matrix, 
which we simply write $I$ if the dimension is clear.
By $\|\cdot\|$ we mean any unitarily invariant matrix norm. 
$\U(n)$ denotes the group of complex unitary matrices.
We let  $\sym X=\frac{1}{2}(X^*+X)$ denote the Hermitian part of $X$ and $\skw X=\frac{1}{2}(X-X^*)$ the skew-Hermitian part of $X$ such that $X=\sym X+\skw X$. $\exp$ denotes the matrix exponential function $\exp X=\sum_{n=0}^\infty \frac{1}{n!}X^n$.
In general, $\Log Z$ with capital letter denotes any solution to $\exp X=Z$, while $\log Z$ denotes the principal matrix logarithm.

\section{Preliminaries}
\subsection{Compound matrices and the generalized Bernstein inequality}
The most important ingredient 
for our proof is inequality \eqref{eq:Cohen} below, which is stated in terms of compound matrices.
The $k$-th compound matrix $A\comp$ of a matrix $A$ is the ${n\choose k}\times{n\choose k}$-matrix consisting of the (lexicographically ordered) determinants of all $k\times k$ submatrices of $A$ (the minors).
For the convenience of the reader we recall some properties of compound matrices (see e.g. \cite[p.411]{Bernstein2009}): 
$$(AB)\comp=A\comp B\comp \qquad \text{for any } A, B\in\Cnn\qquad \text{(Binet-Cauchy formula).}$$
In particular: if $A$ is invertible, $\cdot\comp$ and $\inv$ commute: $$(A\comp)\inv=(A\inv)\comp.$$
Denote by $\trik A:=\tr_i [A\comp]$ the $i$-th partial trace (sum of the $i$ largest eigenvalues 
in modulus) of the $k$-th compound matrix of $A$. If $A$ is similar to $B$, that is $A=SBS\inv$, then 

\begin{align}\label{eq:triksimilar}
\trik A=\trik B,
\end{align}
 because $A\comp$ and $B\comp$ are also similar by the preceding two properties.

For $A=\diag(x_1,\ldots x_n)$, the $k$-th compound matrix $A\comp$ is a diagonal matrix with the different products of $k$ factors $x_i$ as entries.
\begin{ex}
 Let $X=\diag(x_1,x_2,x_3,x_4)$, where $x_1\geq x_2 \geq x_3\geq x_4>0$. Then 
 \begin{align*}
 X^{(1)}&=X,\\X^{(2)}&=\diag(x_1x_2,x_1x_3,x_1x_4,x_2x_3,x_2x_4,x_3x_4),\\ X^{(3)}&=\diag(x_1x_2x_3,x_1x_2x_4,x_1x_3x_4,x_2x_3x_4),\\ X^{(4)}&=\matr{x_1x_2x_3x_4}
 \end{align*}
 and e.g. 
 \begin{align*}
 &\tr X=\tr X^{(1)}=\tr_4^{(1)} X,\\
 &\tr_1^{(1)}X=x_1,\quad \tr_2^{(1)}X=x_1+x_2,\quad \tr_3^{(1)}X=x_1+x_2+x_3,\quad \tr_4^{(1)}X=x_1+x_2+x_3+x_4,\\
 &\tr_1^{(2)}X=x_1x_2,\quad \tr_2^{(2)}X=x_1x_2+x_1x_3,\quad \\
 &\tr_3^{(2)}X=x_1x_2+x_1x_3+\max\{x_1x_4,x_2x_3\},\\
 &\tr_2^{(3)}X=x_1x_2x_3+x_1x_2x_4,\\
 &\trok[1]X=x_1,\quad \trok[2]X=x_1x_2,\quad \trok[3]X=x_1x_2x_3,\quad\trok[4]=x_1x_2x_3x_4.  
\end{align*}
 In general, for 
diagonal matrices $\diag(x_i)$ with $x_1\geq x_2\geq\cdots\geq x_n>0$, one obtains $\trok \diag(x_i)= x_1\cdots x_k$.
\end{ex}

Cohen \cite{Cohen88}, generalizing Bernstein's result \cite{Bernstein88}, 
proved the inequality
\begin{align}\label{eq:Cohen}
 \trik (\exp(A)\exp(A^*))\leq \trik (\exp(A+A^*))
\end{align}
for any $A\in \Cnn$, $k=1,\ldots,n$, $i=1,\ldots,{n\choose k}$ and the matrix-exponential function: $\exp(A)=\sum_{n=0}^\infty \frac{A^n}{n!}$. Of course \eqref{eq:Cohen} is an equality if $AA^*=A^*A$.

We will use the case $i=1$ of these inequalities for compound matrix traces to show the {\bf majorization} of suitable vectors. 
\subsection{Majorization}
Let $x,y\in\Rn$. Then $x$ is said to be {\bf majorized} by $y$, $x \prec y$, if 
\begin{align*}
&\sum_{i=1}^k x_i\sort\leq \sum_{i=1}^k y_i\sort \quad \text{for all } k=1,\ldots, n\\
\text{and }\quad &\sum_{i=1}^n x_i= \sum_{i=1}^n y_i,
\end{align*}
where $x\sort$ denotes the vector $x$ with decreasingly rearranged components.\\
If the latter condition is dropped, we say $x$ is {\bf weakly majorized} by $y$, denoted by $x\prec_w y$, see \cite{Marshall2011inequalities}.\\
\begin{Theorem}(\cite{Marshall2011inequalities})
\label{thm:covexpreservesmajorization}
If $x\prec y$ and $f\colon \R\to\R$ is a convex function, then $(f(x_1),\ldots,f(x_n))\prec_w (f(y_1),\ldots,f(y_n))$. \cite[5.A.1]{Marshall2011inequalities}
\end{Theorem}
This theorem can be proved (see \cite[eqn.~(1.9)]{Ando199417}) by using a characterization of majorization, given in \cite[Thm.8]{HardyLittlewoodPolya29}, via the existence of a doubly stochastic matrix $P$ such that $x=Py$.
We note that the theorem includes Karamata's inequality \cite{Karamata32}, which states $\sum_{i=1}^n f(x_i)\leq \sum_{i=1}^n f(y_i)$ under the same conditions.

Based on an observation of von Neumann \cite{vonNeumann37} (see also \cite[Thm. IV.2.1]{Bhatia97}, \cite[sec.~3.5]{Horn91}) on the relationship between unitarily invariant norms and symmetric gauge functions (norms that are invariant under change of order or signs of components) of their singular values 
and Ky Fan's theorem on a conditition for inequalities of symmetric gauge functions \cite[Thm.4]{Fan51} (see also \cite{bhatia94}),
one has the following important 
connection between majorization and unitarily invariant norms :

\begin{Theorem}(\cite{Fan51})
\label{thm:majorizationandnorms}
Let $X,Y\in\Cnn$ 
be two matrices. Then $$\norm{X}\geq \norm{Y}$$ for all unitarily invariant norms $\norm{\cdot}$ \emph{if and only if} the vectors $\sigma(X), \sigma(Y)$ of singular values 
satisfy $$\sigma(X)\succ_w \sigma(Y).$$
\end{Theorem}

\subsection{Matrix logarithm}\label{sec:matlog}
For every nonsingular $Z\in\GL(n,\C)$ there exists a solution $X\in\C^{n\times n}$ to $\exp X=Z$, which we call a logarithm $X=\Log(\A)$ of $Z$. 
By definition, 
\begin{align*}
 \forall\, X\in \C^{n\times n}:\quad \exp\Log X&=X\,,
\end{align*}
whereas the converse does not have to be true without further assumptions,
\begin{align*}
 \Log\exp X&\neq X\, ,  
\end{align*}
because, as in the scalar case, the matrix logarithm is multivalued depending on the unwinding number~\cite[p.~270]{Higham2008} \cite{unwindingmatrix}: a nonsingular real or complex matrix may have an infinite number of real or complex logarithms. 

If we want to work with one special logarithm with certain desirable properties, 
we use the 
{\bf principal matrix logarithm} $\log X$: 
Let $X\in \C^{n\times n}$, and assume that $X$ has no eigenvalues on $(-\infty,0]$. The principal matrix logarithm of $X$ is the unique logarithm of $X$ (the unique solution $Y\in\C^{n\times n}$ of $\exp Y=X$) whose eigenvalues lie in the strip $\{z\in \C:\; -\pi< \Im(z) <\pi\}$. If $X\in\R^{n\times n}$ and $X$ has no eigenvalues in $(-\infty,0]$, then the principal matrix logarithm is real.

The following statements apply strictly only to the principal matrix logarithm \cite[p.721]{Bernstein2009}:
\begin{align}
\label{rules_principal_logarithm_1}
\log\exp X&=X\quad 
\text{if and only if 
$\abs{\Im (\lambda)}<\pi$ for all $\lambda\in \rm{spec}(X)$}
\, , \\
 \log(Q^* X Q)&=Q^* \log(X) \, Q\, ,\quad \forall\, Q\in\U(n)\, .\label{rules_principal_logarithm_2}
\end{align}
Since $\sym X$ is Hermitian the matrix $\exp\sym X$ is positive definite, so we can apply $\eqref{rules_principal_logarithm_1}$ and it follows from $\eqref{rules_principal_logarithm_2}$ that
  \begin{align}    
  \label{log_on_psym_2}
 \forall\,X\in \C^{n\times n}:\quad Q^*[\sym X]Q= Q^*[\log\exp\sym X] Q&=\log[Q^*(\exp\sym X) Q]  \, .
\end{align}

\section{The minimization}\label{sec:main}
\subsection{Preparation}
The goal is to find the unitary $Q\in\U(n)$ that minimizes $\|\Log(Q^* \A)\|$ and $\norm{\sym\Log (Q^* \A)}$ over all possible logarithms.
Due to the non-uniqueness of the logarithm, we give the following as the statement of the minimization problem: 
\begin{align}
\label{min_defi}
\min_{Q\in\U(n)}\norm{\Log (Q^* \A)}&:=\min_{Q\in\U(n)}\{  \norm{X}\in\R\,|\, \exp X=Q^* \A\}\, ,\notag\\
   \min_{Q\in\U(n)}\norm{\sym\Log (Q^* \A)}&:=\min_{Q\in\U(n)}\{  \norm{\sym X}\in\R\,|\, \exp X=Q^* \A\}\, .
\end{align}
We first observe, as shown in \cite{NNF}, 
 that without loss of generality we may assume that $\A\in \GL(n,\C)$ is real, diagonal and positive definite. To see this, consider the unique polar decomposition $\A=U_p \, H$ 
and the eigenvalue decomposition 
$H=V D V^*$ 
where $D=\diag(d_1,\ldots, d_n)$ with $d_i> 0$. 
Then
\begin{align}
 \min_{Q\in\U(n)}\norm{\sym\Log (Q^* \A)}&=\min_{Q\in\U(n)}\{  \norm{\sym X} \,|\, \exp X=Q^* \A\}\notag\\
&=\min_{Q\in\U(n)}\{  \norm{\sym X} \,|\, \exp X=Q^*  U_pV D V^*\}\notag\\
&=\min_{Q\in\U(n)}\{  \norm{\sym X} \,|\, \exp (V^*XV)=V^* Q^* U_pV D\}\notag\\
&=\min_{\widetilde{Q}\in\U(n)}\{  \norm{\sym X} \,|\,  \exp (V^*XV)=\widetilde{Q}^* D \}   \\
&=\min_{\widetilde{Q}\in\U(n)}\{  \norm{\sym (V^*X V)} \,|\,  \exp (V^*XV)=\widetilde{Q}^* D \}\notag\\
&=\min_{\widetilde{Q}\in\U(n)}\{  \norm{\sym (\widetilde{X}) } \,|\,  \exp (\widetilde{X})=\widetilde{Q}^* D \}\notag\\
&= \min_{Q\in\U(n)}\norm{\sym\Log Q^*D}\, ,\notag
 \end{align}
where we used the unitary invariance for any unitarily invariant matrix norm and the fact that $X\mapsto \sym X$ and $X\mapsto \exp X$ are isotropic matrix functions, i.e.  $f(V^* X V)=V^* f(X) V$ for all unitary $V$. If the minimum is achieved for $Q=\id$ in $\min_{Q\in \U(n)}\norm{\sym \Log (Q^* D) }$ then this corresponds to $Q=U_p$ in $\min_{Q\in \U(n)}\norm{\sym \Log Q^* \A}$. Therefore, in the following we assume  that $Z=D=\diag(d_1,\ldots d_n)$ with $d_1\ge d_2\ge \ldots \ge d_n>0$.

\subsection{Main result: minimizing $\|$Log Q$^*$D$\|$}
\label{sec:minimizingHERE}
Our starting point is the problem of minimizing the Hermitian part
\begin{align*}
\min_{Q\in\U(n)}\|\sym(\Log(Q^* \A))\|. 
\end{align*}
 As we will see, a solution of this problem will already imply the other minimization properties.
Let $n\in\N$ be arbitrary. 
For any $Q\in \U(n)$ the Hermitian positive definite matrix $\exp(\sym\Log Q^* D)$ can be unitarily diagonalized with positive, real eigenvalues, i.e., for some $Q_1\in\U(n)$
\begin{align}
  \label{eigenvalue_representation}
    Q_1^* \left(\exp(\sym\Log Q^* D)\right) Q_1
    =\diag(x_1,\ldots,x_n)=:X   \, .
\end{align}
Here, we assume that the positive real eigenvalues are ordered as  $x_1\ge x_2\ge\ldots \ge x_n>0$. 
Then from 
\begin{align*}
\det\exp X &=e^{\tr X}\, , 
\end{align*}
which holds true for arbitrary matrices \text{\cite[p.712]{Bernstein2009}}, we have
\begin{align}
 \det X &= \det(\exp(\sym\Log Q^*D))=|\exp(\tr(\sym\Log Q^*D))|\notag\\
 &=\exp(\Re\tr(\Log Q^*D))=|\exp(\tr(\Log Q^*D))|\\
 &=|\det \exp\Log(Q^*D)|=|\det(Q^*D)|=\det D,\notag
\end{align}
and therefore
\begin{align}
\label{eq:detequal}
 x_1\,x_2\cdots x_{n-1}\,x_n= d_1\,d_2\cdots d_{n-1}\,d_n.
\end{align}

Due to \eqref{eigenvalue_representation}, $X^2=XX^*$ and $\exp(\sym\Log Q^*D) \exp(\sym\Log Q^*D)^*$ are similar. Furthermore,
\begin{align*} 
\exp(\sym\Log Q^*D) \exp(\sym\Log Q^*D)^*&=\exp(\sym (\Log Q^*D))\exp((\sym(\Log Q^*D))^*)\\
&=\exp(2\sym\Log Q^*D) \\&=\exp((\Log  Q^*D)+(\Log(Q^*D))^*).
\end{align*}
Hence by equation \eqref{eq:triksimilar} the matrices $X^2$ and $\exp((\Log  Q^*D)+(\Log(Q^*D))^*)$ have the same partial compound traces $\trik$, and setting $i=1$ we obtain from Cohen's inequality \eqref{eq:Cohen} 
that 
\begin{align*}
 \trok X^2 &= \trok (\exp((\Log  Q^*D)+(\Log(Q^*D))^*))\\
 &\geq \trok (\exp(\Log Q^*D)\exp(\Log Q^*D)^*)=\trok(Q^*DD^*Q)=\trok D^2,
\end{align*}
i.e. (recall that $\trok$ is the largest eigenvalue of the $k$-th compound matrix):
\begin{align}
\label{eq:ximajdi}
 x_1^2&\geq d_1^2\notag\\
 x_1^2\,x_2^2&\geq d_1^2\,d_2^2\notag\\
 &\;\,\vdots\\
 x_1^2\,x_2^2\cdots x_{n-1}^2&\geq d_1^2\,d_2^2\cdots d_{n-1}^2\notag\\
 x_1^2\,x_2^2\cdots x_{n-1}^2\,x_n^2&\geq d_1^2\,d_2^2\cdots d_{n-1}^2\,d_n^2.\notag
\end{align}
Of course, by \eqref{eq:detequal} the last inequality is in fact an equality.
Applying the logarithm to \eqref{eq:detequal} and to \eqref{eq:ximajdi} gives
\begin{align}
 \log x_1&\geq \log d_1\notag \\
 \log x_1+\log x_2&\geq \log d_1+ \log d_2\notag \\
 &\;\,\vdots\notag \\
 \log x_1+\log x_2+\cdots +\log x_{n-1}&\geq \log d_1+\log d_2+\cdots +\log d_{n-1}\label{logmajor}\\
 \log x_1+\log x_2+\cdots +\log x_{n}&= \log d_1+\log d_2+\cdots +\log d_{n}.\notag 
\end{align}
That is, we have the majorization 
\begin{align}\label{logmajor}
(\log x_1,\ldots,\log x_n)\succ (\log d_1,\ldots,\log d_n), \end{align}
very much in the spirit of the reformulation of Cohen's result in \cite[Thm. C]{AndoHiai}.\\
As the modulus is a convex function, from Theorem \ref{thm:covexpreservesmajorization}  we obtain 
\begin{align}\label{abslogmajor}
(|\log x_1|,\ldots,|\log x_n|)\succ_w (|\log d_1|,\ldots,|\log d_n|)\,.
 \end{align}
Note that these vectors contain nothing but the singular values of $\log X$ and $\log D$ respectively, and hence 
\begin{align}
\label{eq:norminequality}
 \norm{\log X}\geq \norm{\log D} 
\end{align}
for any unitarily invariant norm, by Theorem \ref{thm:majorizationandnorms}.

According to \eqref{log_on_psym_2} and \eqref{eigenvalue_representation},
\begin{align}
 \log X = \log[Q_1^* \exp(\sym\Log Q^* D) Q_1]=Q_1^* (\sym\Log Q^* D) Q_1
\end{align}
and because $\norm\cdot$ is unitarily invariant, 
\eqref{eq:norminequality} can be stated as 
\begin{align}
 \norm{\sym \Log Q^*D}\geq \norm{\log D}.
\end{align}
Together with the trivial upper bound (let $Q=\id$ and $\Log=\log$) we conclude that
\begin{align}\label{symlogmin}
\min_{Q\in \U(n))}\norm{\sym \Log (Q^* D) }=\norm{\log D}\, .
\end{align}
The minimum is realized for $Q=\id$, which corresponds to the polar factor $U_p$ in the original formulation.

To obtain the solution for the minimization problem $\min_{Q\in \U(n))}\norm{\Log (Q^* D) }$ from that of the Hermitian part \eqref{symlogmin}, we use the fact that
for any unitarily invariant norm, the norm of the Hermitian part of any matrix is less than or equal to the norm of the matrix \cite[p.454]{Horn85}, 
\[
 \norm{\sym X}\leq \norm{X}. 
\]
It follows that
\begin{align*}
 \min_{Q\in \U(n)}\norm{\Log Q^*D}\geq\min_{Q\in \U(n)}\norm{\sym \Log Q^*D}=\norm{\log D}\,.
\end{align*}
The last inequality, together with the upper bound for $Q=\id$, yields
\begin{equation}  \label{eq:minprop}
 \min_{Q\in \U(n)}\norm{\Log (Q^* D) }=\norm{\log D}\, .  
\end{equation}

Combining the above results, 
for all $\my>0, \my_c\geq 0$ we obtain \eqref{eq:minimization_family}:
\begin{align*}
 \min_{Q\in \U(n)}& \my\norm{\sym\Log(Q^*Z)}+\my_c\norm{\skw\Log(Q^*Z)} \geq \min_{Q\in \U(n)} \my\norm{\sym\Log(Q^*Z)} \\
 &=\my\norm{\sym\Log(U_p^*Z)}=\my\norm{\sym\Log(U_p^*Z)}+\my_c\underbrace{\norm{\skw \Log(U_p^*Z)}}_{=0}\\
 &=\my\norm{\log U_p^*Z}=\my\norm{\log\sqrt{Z^*Z}}\,. 
\end{align*}
In summary, we have proved the following:
\begin{Theorem}\label{thm:main}
 Let $Z\in \Cnn$ be a nonsingular matrix and let $Z=U_pH$ be its polar decomposition. Then for any unitarily invariant norm $\|\cdot\|$
\[
 \min_{Q\in U(n)} \norm{\Log Q^*Z}=\norm{\log U_p^*Z}=\norm{\log H},
\]
\[
 \min_{Q\in U(n)} \norm{\sym\Log Q^*Z}=\norm{\sym\log U_p^*Z}=\norm{\log H}, \]
and for any $ \my>0,\my_c\geq 0$
\[ 
\min_{Q\in U(n)}\left(\my\norm{\sym\Log(Q^*Z)}+\my_c\norm{\skw\Log(Q^*Z)}\right)
= \mu\norm{\log H}
.\]
\end{Theorem}
\subsection{Uniqueness}\label{sec:unique}
The question of uniqueness was considered in \cite{NNF} for the spectral norm and for the Frobenius norm when $n\leq 3$. 
The analysis there showed that 
$Q=U_p$ is the unique minimizer of $\|\Log Q^*Z\|$ for the Frobenius norm for $n\leq 3$, but not for the spectral norm. 
Moreover, it was conjectured there that $Q=U_p$ is the \emph{only} matrix that minimizes $\|\log Q^*Z\|$ 
regardless of the choice of the unitarily invariant norm. 
Here we prove this in the affirmative:
\begin{Theorem}
 Let $Z\in \Cnn$ be a nonsingular matrix, and 
 suppose $\widehat{Q}\in U(n)$ is such that for every unitarily invariant norm $\norm{\cdot}$ the equality
 \[
  \norm{\log \widehat{Q}^*Z}=\min_{Q\in U(n)}\norm{\Log Q^*Z}
 \]
 holds. Then $\widehat{Q}=U_p$, the unitary polar factor of $Z$.
\end{Theorem}
{\sc proof}. 
By Theorem~\ref{thm:majorizationandnorms}, for a fixed $Q$ to be the minimizer of 
$\|$Log Q$^*$Z$\|$ for every unitarily invariant norm, 
we need equality to hold in \eqref{eq:norminequality} 
for every Ky Fan $k$-norm $\|Z\|_{(k)}=\sum_{i=1}^k\sigma_i(Z)$ for $k=1,2,\ldots,n$. 
That is, we require
\begin{equation}  \label{eq:kyfans}
 \norm{\log X}_{(k)}= \norm{\log D}_{(k)},\quad k=   1,2,\ldots,n. 
\end{equation}

We re-order the sets $\log x_i$ and $\log d_i$ to arrange in 
decreasing order of absolute value and denote them by
$|\log \widehat x_1|\geq |\log \widehat x_2|\geq \cdots \geq |\log \widehat x_n|$ and 
$|\log \widehat d_1|\geq |\log \widehat d_2|\geq \cdots \geq |\log \widehat d_n|$. 
Then \eqref{eq:kyfans}
is equivalent to 
\begin{equation}  \label{eq:eqmaj}
(|\log \widehat x_1|,\ldots,|\log \widehat x_n|)= (|\log \widehat d_1|,\ldots,|\log \widehat d_n|).   
\end{equation}
Recall that $\log x_i$ and $\log d_i$ also satisfy the majorization property \eqref{logmajor}. 
We now claim that \eqref{logmajor} and \eqref{eq:eqmaj} imply $x_i=d_i$ for all $i=1,\ldots, n$. 

It is worth noting that \eqref{eq:eqmaj} includes the statement \[\sum_{i=1}^n|\log \widehat x_i|=\sum_{i=1}^n|\log \widehat d_i|\Leftrightarrow \sum_{i=1}^n|\log x_i|=\sum_{i=1}^n|\log d_i|\,,\] that is, Karamata's inequality holds with equality. 
Moreover, Karamata's inequality is known to become an equality if and only if the two sets are equal, which in this case means
$\log x_i=\log d_i$ for all $i$, provided that the function $f(x)$ (which here is $|x|$) is strictly convex. 
However, since $|x|$ is not strictly convex over $\R$, this argument is not directly applicable. Below we shall see that we nonetheless have $x_i=d_i$.

First, since $\log$ is a monotone function we have either $\log\widehat x_1=\log x_1\geq 0$ or $\log \widehat x_1=\log x_n<0$, and similarly $\log \widehat d_1=\log d_1\geq 0$ or $\log \widehat d_1=\log d_n<0$. 

By \eqref{eq:eqmaj} we need $|\log\widehat x_1| = |\log\widehat d_1|$, so 
either $\log\widehat x_1 = \log\widehat d_1$ or 
$\log\widehat x_1 = -\log\widehat d_1$. 
Now if $\widehat x_1 = \widehat d_1$ then we can remove $\widehat x_1,\widehat d_1$ from the lists $\widehat x_i,\widehat d_i$
without affecting the argument. 
Hence here we suppose that 
$\log\widehat x_1 = -\log\widehat d_1$, and show by contradiction that this assumption cannot hold, thus proving 
$\widehat x_1 = \widehat d_1$.

Observe that the assumption $\log\widehat x_1 = -\log\widehat d_1$ forces $\widehat x_1=x_1$ and $\widehat d_1=d_n$
(instead of $\widehat x_1=x_n$ or $\widehat d_1=d_1$), 
because if $\log \widehat x_1=\log x_n<0$ and hence $|\log \widehat x_1|>|\log x_1|$, 
then $\log\widehat d_1>0$ and hence $|\log\widehat d_1|=\log d_1\geq 0$, 
and so $\log d_1=|\log\widehat d_1|=|\log\widehat x_1|>|\log x_1|$, contradicting the first majorization property in \eqref{logmajor}. 

Hence our assumptions are $\log\widehat x_1 = -\log\widehat d_1,\,\,\widehat x_1=x_1$ and $\widehat d_1=d_n$. 
In the  equality $\sum_{i=1}^n\log x_i=\sum_{i=1}^n\log d_i$ of \eqref{logmajor}, 
 subtracting $\log d_{n}= \log \widehat d_1=-\log x_1$ from both sides yields
\[
 2\log x_1 + \log x_2 +\cdots + \log x_{n} =  \log d_1 + \log d_2 +\cdots + \log d_{n-1}. 
\]
Together with the $(n-1)$th majorization assumption $\sum_{i=1}^{n-1}\log x_i\geq \sum_{i=1}^{n-1}\log d_i$ we need 
\[
 \log x_1 + \log x_2 +\cdots + \log x_{n-1} \geq  
2\log x_1 + \log x_2 +\cdots + \log x_{n},
\]
which is equivalent
to $ \log x_1+\log x_n\leq 0$. This contradicts our assumption $|\log x_1|\geq|\log x_n|$ unless
$|\log x_1|=|\log x_n|$, but in this case we can remove both $x_n$ and $d_n$ (with $x_n=d_n$) from the list without affecting the argument. Overall we have shown that 
 we have $x_1=d_1$, and by repeating the same argument we conclude that 
\begin{align}\label{eq:xdequality} x_i=d_i \qquad\text{for all }i. \end{align}

%


We next  examine the necessary conditions to satisfy $\norm{\sym \Log (Q^* D) }=\norm{\log D}$ in \eqref{symlogmin}, and show that we need $Q=I$. 
 We clearly need \[\norm{\sym \Log (Q^* D) }=\norm{\Log (Q^* D) }\] for every unitarily invariant norm, which forces $\Log (Q^* D)$ to be Hermitian. 
Hence the matrix $\exp(\Log (Q^* D))$ is positive definite, 
so we can write $\exp(\Log (Q^* D)) = Q_1^*\diag(x_1,\ldots,x_n)Q_1$ 
 for some unitary $Q_1$ and $x_i>0$. Therefore 
the matrix logarithm is necessarily the principal one, and 
 \begin{align}\label{eq:qd}
\Log (Q^* D) =  \log (Q^* D) = Q_1^*\diag(\log x_1,\ldots,\log x_n)Q_1.    
 \end{align}
Hence by \eqref{eq:qd} and \eqref{eq:xdequality} we have \[\log (Q^* D) = Q_1^*\diag(\log x_1,\ldots,\log x_n)Q_1 = Q_1^*\log(D)Q_1,\] so 
taking the exponential of both sides yields
\begin{align}\label{eq:qd2}
Q^* D = Q_1^*DQ_1.  
\end{align}
Hence $D = Q(Q_1^*DQ_1).$ Note that this is the polar decomposition of $D$, as $Q_1^*DQ_1$ is Hermitian positive definite. It follows that $Q$ must be equal to the unique unitary polar factor of $D$, which is clearly $I$. 
Overall,  for \eqref{eq:qd2} to hold we always need $Q = I$, which corresponds to the unitary polar factor $U_p$ in the original formulation. 
Thus $Q = U_p$ is the unique minimizer of $\norm{\Log (Q^* D) }$ with minimum $\norm{\log (U_p^* D) }$. Other choices of the matrix logarithm are easily seen to give larger $\norm{\Log (Q^* D) }$. 
\hfill$\square$
\medskip

Although we have shown that $Q=U_p$ is always a minimizer of $\|\Log Q^*Z\|$, for a specific unitarily invariant norm it may not be the unique minimizer. For example, for the spectral norm there can be infinitely many $Q$ for which $\|\log Q^*Z\|=\|\log U_p^*Z\|$, as  was shown in \cite{NNF}. 
In general, $Z=U_p$ is not the unique minimizer when the norm does not involve all the singular values, such as the spectral norm $\|Z\|=\sigma_1(Z)$ and Ky Fan $k$-norm $\|Z\|=\sum_{i=1}^k\sigma_i(Z)$ for $k< n$. 

Below we discuss a general form of the minimizers $Q$ for a Ky fan $k$-norm. 
  \begin{Proposition}
$\A=U\Sigma V^*$ be an SVD of a nonsingular $\A$ with 
$\Sigma = \diag(\sigma_1,\sigma_2,\ldots, \sigma_n)$. 
Let $\widehat\Sigma = \diag(\widehat\sigma_1,\widehat\sigma_2,\ldots, \widehat\sigma_n)$ where $\{\widehat\sigma_i\}$ is a permutation of $\{\sigma_i\}$ such that 
$|\log\widehat\sigma_1|\geq  \cdots\geq |\log\widehat\sigma_n|$, and define $\widehat U,\widehat V$ such that $\A=\widehat U\widehat\Sigma \widehat V^*$ is an SVD with permuted order of singular values. 


Then for any $\widehat Q\in U(n)$ expressed as
$\widehat Q = \widehat U\diag(I_k,Q_{22})\widehat V^*$ where $Q_{22}\in U(n-k)$ such that 
\begin{equation}  \label{eq:q22cond}
\|\log Q_{22}\diag(\widehat\sigma_{k+1},\ldots, \widehat\sigma_n)\|_{2}\leq |\log\widehat\sigma_{k}|, 
\end{equation}
 (where $\|\cdot\|_{2}$ denotes the spectral norm, that is, the largest singular value), we have
 \begin{equation}   \label{eq:q2goal}
\norm{\log \widehat{Q}^*Z}_{(k)}=\min_{Q\in U(n)} \norm{\Log Q^*Z}_{(k)}.   
 \end{equation}
  \end{Proposition}
{\sc proof}. 
Direct calculation shows for such $\widehat Q$ that 
\[\log \widehat Q^*\A =  \log \widehat V \diag(I_k,Q_{22} )\widehat \Sigma \widehat V^*, \]
so the singular values of $\log \widehat Q^*\A$ are the union of $\log\widehat\sigma_i$, $i=1,\ldots, k$ and 
those of $\log \left(Q_{22}\diag(\widehat\sigma_{k+1},\ldots, \widehat\sigma_n)\right)$. 
By \eqref{eq:q22cond} we have $\norm{\log \widehat{Q}^*Z}_{(k)}=\sum_{i=1}^k\log\widehat\sigma_i=\norm{\log U_p^*Z}_{(k)}$, and 
\eqref{eq:q2goal} follows from the fact that \[\norm{\log U_p^*Z}_{(k)}=\min_{Q\in U(n)} \norm{\Log Q^*Z}_{(k)}\] as we have seen in Theorem~\ref{thm:main}. 
\hfill$\square$

\medskip

We note that the set of $Q_{22}$ that satisfies \eqref{eq:q22cond}  includes the choice $Q_{22}=I_{n-k}$. 
Moreover, the set generally includes more than $I_{n-k}$, and can be (but not always) as large as the whole group $U(n-k)$.

%
%
%
%

\subsection{Rectangular $Z$}
The polar decomposition $Z=U_p\,H$ is defined for any $Z\in\Cmn$ with $m\geq n$, including singular and rectangular matrices~\cite[Ch.~8]{Higham2008}. Also in this case it solves \cite[Thm.~8.4]{Higham2008} 
\[
 \norm{Z-U_p}=\min\{\norm{Z-Q}: Q^*Q=I_n\}.
\]
Therefore a natural question arises of whether $U_p$ is still the minimizer of 
$\|\Log Q^*Z\|$ over $Q\in\Cmn$ such that $Q^*Q=I_n$ when $m>n$. 

The answer to this question is in the negative, as can be seen by the simple example 
$Z=\big[\begin{smallmatrix}
1\cr 1  \end{smallmatrix}\big]$, 
for which $Z = U_pH$ with 
$U_p=\frac{1}{\sqrt{2}}\big[\begin{smallmatrix}1\cr 1  \end{smallmatrix}\big]$ and $H=\sqrt{2}$. 
Defining $V  =\big[\begin{smallmatrix}1\cr 0  \end{smallmatrix}\big]$
we have $\log U_p^*Z=\frac{1}{\sqrt{2}}$ but 
$\log V^*Z=0$, clearly showing that $U_p$ is generally not the minimizer of $\|\Log Q^*Z\|$. 
We conclude that the minimization property of $U_p$ that we have discussed is particular for square and nonsingular matrices, contrary to the minimization property of $U_p$ with respect to $\norm{Z-Q}$, which holds for any $Z$ including rectangular ones.

{\footnotesize
    \bibliographystyle{plain} 

    \bibliography{literatur1}

\begin{thebibliography}{10}

\bibitem{Ando199417}
T.~Ando.
\newblock Majorizations and inequalities in matrix theory.
\newblock {\em Linear Algebra Appl.}, 199, Supplement 1(0):17 -- 67, 1994.
\newblock Special issue honoring Ingram Olkin.

\bibitem{AndoHiai}
T.~Ando and F.~Hiai.
\newblock Log majorization and complementary {G}olden-{T}hompson type
  inequalities.
\newblock {\em Linear Algebra Appl.}, 197-198(0):113 -- 131, 1994.

\bibitem{unwindingmatrix}
M.~Aprahamian and N.~J. Higham.
\newblock The matrix unwinding function, with an application to computing the
  matrix exponential.
\newblock {MIMS EPrint} 2013.21, The University of Manchester, UK, May 2013.

\bibitem{Autonne1902}
L.~Autonne.
\newblock Sur les groupes lin{\'e}aires, r{\'e}els et orthogonaux.
\newblock {\em Bull. Soc. Math. France}, 30:121--134, 1902.

\bibitem{Bernstein88}
D.~S. Bernstein.
\newblock Inequalities for the trace of matrix exponentials.
\newblock {\em SIAM J. Matrix Anal. Appl.}, 9(2):156--158, 1988.

\bibitem{Bernstein2009}
D.~S. Bernstein.
\newblock {\em Matrix Mathematics}.
\newblock Princeton University Press, New Jersey, 2009.

\bibitem{Bhatia97}
R.~Bhatia.
\newblock {\em Matrix Analysis}, volume 169 of {\em Graduate Texts in
  Mathematics}.
\newblock Springer, New-York, 1997.

\bibitem{bhatia94}
R.~Bhatia and C.~Davis.
\newblock Relations of linking and duality between symmetric gauge functions.
\newblock In A.~Feintuch and I.~Gohberg, editors, {\em Nonselfadjoint Operators
  and Related Topics}, volume~73 of {\em Operator Theory: Advances and
  Applications}, pages 127--137. Birkhaeuser Basel, 1994.

\bibitem{Bhatia2012726}
R.~Bhatia and P.~Grover.
\newblock Norm inequalities related to the matrix geometric mean.
\newblock {\em Linear Algebra Appl.}, 437(2):726 -- 733, 2012.

\bibitem{NeffBirsan13}
M.~B\^{i}rsan and P.~Neff.
\newblock Existence of minimizers in the geometrically non-linear 6-parameter
  resultant shell theory with drilling rotations.
\newblock {\em Mathematics and Mechanics of Solids}, 2013.

\bibitem{Neff_log_inequality13}
M.~B\^{i}rsan, P.~Neff, and J.~Lankeit.
\newblock Sum of squared logarithms: {A}n inequality relating positive definite
  matrices and their matrix logarithm.
\newblock {\em Journal of Inequalities and Applications}, 2013(1):168, 2013.

\bibitem{ByXu08}
R.~Byers and H.~Xu.
\newblock {A new scaling for Newton's iteration for the polar decomposition and
  its backward stability}.
\newblock {\em SIAM J. Matrix Anal. Appl.}, {30}:{822--843}, {2008}.

\bibitem{Cohen88}
J.~E. Cohen.
\newblock Spectral inequalities for matrix exponentials.
\newblock {\em Linear Algebra Appl.}, 111:25--28, 1988.

\bibitem{Fan51}
K.~Fan.
\newblock Maximum properties and inequalities for the eigenvalues of completely
  continuous operators.
\newblock {\em Proc. Natl. Acad. Sci. USA.}, 37(11):760 -- 766, 1951.

\bibitem{Fan55}
K.~Fan and A.~J. Hoffmann.
\newblock Some metric inequalities in the space of matrices.
\newblock {\em Proc. Amer. Math. Soc.}, 6:111--116, 1955.

\bibitem{Golubbookori}
G.~H. Golub and C.~V.~Van Loan.
\newblock {\em Matrix Computations}.
\newblock {The Johns Hopkins University Press}, 1996.

\bibitem{HardyLittlewoodPolya29}
G.H. Hardy, J.E. Littlewood, and G.~P\'{o}lya.
\newblock Some simple inequalities satiesfied by convex functions.
\newblock {\em Messenger of Mathematics}, 58:145 -- 152, 1929.

\bibitem{Higham2008}
N.~J. Higham.
\newblock {\em Functions of {M}atrices: Theory and {C}omputation}.
\newblock SIAM, Philadelphia, PA, USA, 2008.

\bibitem{Horn85}
R.~A. Horn and C.~R. Johnson.
\newblock {\em Matrix {A}nalysis.}
\newblock Cambridge University Press, New York, 1985.

\bibitem{Horn91}
R.~A. Horn and C.~R. Johnson.
\newblock {\em Topics in {M}atrix {A}nalysis.}
\newblock Cambridge University Press, New York, 1991.

\bibitem{Karamata32}
J.~Karamata.
\newblock Sur une in\'{e}galit\'{e} relative aux fonctions convexes.
\newblock {\em Publ. Math. Univ. Belgrad}, 1:145--148, 1932.

\bibitem{Marshall2011inequalities}
A.~W. Marshall, I.~Olkin, and B.~C. Arnold.
\newblock {\em Inequalities: Theory of Majorization and Its Applications}.
\newblock Springer Series in Statistics. Springer, 2011.

\bibitem{nakatsukasa:2700}
Y.~Nakatsukasa, Z.~Bai, and F.~Gygi.
\newblock {Optimizing Halley's iteration for computing the matrix polar
  decomposition}.
\newblock {\em SIAM J. Matrix Anal. Appl.}, 31(5):2700--2720, 2010.

\bibitem{nahi11}
Y.~Nakatsukasa and N.~J. Higham.
\newblock Backward stability of iterations for computing the polar
  decomposition.
\newblock {\em SIAM J. Matrix Anal. Appl.}, 33(2):460--479, 2012.

\bibitem{nakahig12}
Yuji Nakatsukasa and Nicholas~J. Higham.
\newblock Stable and efficient spectral divide and conquer algorithms for the
  symmetric eigenvalue decomposition and the svd.
\newblock {\em SIAM J. Sci. Comp}, 35(3):A1325--A1349, 2013.

\bibitem{Neff01c}
P.~Neff.
\newblock Local existence and uniqueness for quasistatic finite plasticity with
  grain boundary relaxation.
\newblock {\em Quart. Appl. Math.}, 63:88--116, 2005.

\bibitem{Neff_micromorphic_rse_05}
P.~Neff.
\newblock Existence of minimizers for a finite-strain micromorphic elastic
  solid.
\newblock {\em Proc. Roy. Soc. Edinb. A}, 136:997--1012, 2006.

\bibitem{Neff_Chelminski_ifb07}
P.~Neff and K.~Che{\l}mi\'nski.
\newblock A geometrically exact {C}osserat shell-model for defective elastic
  crystals. {J}ustification via {$\Gamma$}-convergence.
\newblock {\em Interfaces and Free Boundaries}, 9:455--492, 2007.

\bibitem{NeffEidelOsterbrinkMartin_Riemannianapproach}
P.~{Neff}, B.~{Eidel}, F.~{Osterbrink}, and R.~{Martin}.
\newblock {A Riemannian approach to strain measures in nonlinear elasticity}.
\newblock {\em ArXiv e-prints}, May 2013.

\bibitem{Neff_Osterbrink_Martin_hencky13}
P.~Neff, B.~Eidel, F.~Osterbrink, and R.~Martin.
\newblock {T}he isotropic {H}encky strain energy measures the geodesic distance
  of the deformation gradient {$F \in\mathrm{GL^+}(n)$} to $\mathrm{SO}(n)$ in
  the unique left invariant {R}iemannian metric on $\mathrm{GL}(n)$ which is
  also right $\mathrm{O}(n)$-invariant.
\newblock {\em submitted}, 2013.

\bibitem{Neff_Biot07}
P.~Neff, A.~Fischle, and I.~M\"unch.
\newblock Symmetric {C}auchy-stresses do not imply symmetric {B}iot-strains in
  weak formulations of isotropic hyperelasticity with rotational degrees of
  freedom.
\newblock {\em Acta Mechanica}, 197:19--30, 2008.

\bibitem{Neff_Forest_jel05}
P.~Neff and S.~Forest.
\newblock A geometrically exact micromorphic model for elastic metallic foams
  accounting for affine microstructure. {M}odelling, existence of minimizers,
  identification of moduli and computational results.
\newblock {\em J. Elasticity}, 87:239--276, 2007.

\bibitem{NNF}
P.~{Neff}, Y.~{Nakatsukasa}, and A.~{Fischle}.
\newblock The unitary polar factor ${Q}={U}_p$ minimizes $\|
  \rm{{L}og}$(${Q}^*${Z})$\|^2$ and $\|\rm{sym_*\, {L}og}$(${Q}^*${Z})$\|^2$ in
  the spectral norm in any dimension and the {F}robenius matrix norm in three
  dimensions.
\newblock {\em to appear in SIMAX}, 2013.
\newblock http://arxiv.org/abs/1302.3235.

\bibitem{vonNeumann37}
J.~V. Neumann.
\newblock Some matrix-inequalities and metrization of matric-space.
\newblock {\em Tomsk University Review}, 1:286 -- 300, 1937.
\newblock reprinted in A.H. Taub (Ed.), John Von Neumann collected works, Vol.
  4 Pergamon, New York (1962).

\end{thebibliography}
}    
    
\section{Appendix}
\subsection{Additional conditions to be gained from Cohen's formula}
Cohen's inequality \eqref{eq:Cohen} gave us the estimates we needed for the majorization of the singular values of the matrix logarithms - but in fact it gives us also some more inequalities that we did not use.
In \cite{Neff_log_inequality13} (on which \cite{NNF} relies) conditions for the sum of squared logarithms inequality \cite{Neff_log_inequality13}
\begin{align}
\label{eq:sumofsquaredlog}
 \log^2(y_1)+\log^2(y_2)+\log^2(y_3) \geq \log^2(a_1)+\log^2(a_2)+\log^2(a_3)
\end{align}
for $n=3$ to hold are stated in terms of elementary symmetric polynomials (for the case $n=3$)
  \begin{align}
 \label{eq:conditionsforsumofsqaredlog}
 e_0^{(3)}(y_1,y_2,y_3)&=1 &\geq e_0^{(3)}(a_1,a_2,a_3)\;\notag\\
 e_1^{(3)}(y_1,y_2,y_3)&=y_1+y_2+y_3 &\geq e_1^{(3)}(a_1,a_2,a_3)\;\\
 e_2^{(3)}(y_1,y_2,y_3)&=y_1y_2+y_1y_3+y_2y_3 &\geq e_2^{(3)}(a_1,a_2,a_3)\; \notag\\
 e_3^{(3)}(y_1,y_2,y_3)&=y_1y_2y_3 & \geq e_3^{(3)}(a_1,a_2,a_3). \notag 
  \end{align}
The sum of squared logarithms inequality states that \eqref{eq:conditionsforsumofsqaredlog} implies \eqref{eq:sumofsquaredlog}.

We will see that these conditions \eqref{eq:conditionsforsumofsqaredlog} also arise from Cohen's inequality \eqref{eq:Cohen}.
The assertion of \eqref{eq:Cohen} is that $\trik X^2\geq \trik D^2$ (where $X=\diag(x_1,\ldots,x_n)$ is similar to $\exp(\sym\Log Q^* D)$ and $D=\diag(d_1,\ldots,d_n)$, where in both cases the eigenvalues arranged in descending order). 

To save some squares, we let $y_i=x_i^2$ und $a_i=d_i^2$. Then the following expressions are greater than their counterparts with $a$ instead of $y$:
\begin{align}
 y_1&\geq a_1\notag\\
 y_1+y_2 &\geq a_1+a_2\notag\\
 y_1+y_2+y_3&\geq a_1+a_2+a_3\\
 y_1+y_3+y_3+y_4 &\geq a_1+a_2+a_3+a_4\notag\\
 \vdots\notag\\
 e^{(n)}_1(y)=y_1+\cdots+y_n&\geq a_1+\cdots+a_n=e^{(n)}_1(a),\notag
\end{align}
where $e^{(n)}_i(y)=e_i(y_1,\cdots,y_n)$ denotes the elementary symmetric polynomial of first order, and
\begin{align}
 y_1y_2 & \geq a_1a_2\notag\\
 y_1y_2+y_1y_3 &\geq a_1a_2+a_1a_3\\
 y_1y_2+y_1y_3+y_1y_4&\geq a_1a_2+a_1a_3+a_1a_4\,.\notag
 \end{align}
We note that the last inequality $y_1y_2+y_1y_3+y_1y_4\geq a_1a_2+a_1a_3+a_1a_4$ does not necessarily hold in this form. All of the appearing sums have to be the sums of the greatest corresponding terms.
For example, if $y_2y_3$ is greater than $y_1y_4$ (and nothing excludes that), then the left hand side becomes $y_1y_2+y_1y_3+y_2y_3$. 
Whether the right hand side stays $a_1a_2+a_1a_3+a_1a_4$ or is changed to $a_1a_2+a_1a_3+a_2a_3$ also, depends on, whether $a_2a_3$ or $a_1a_4$ is larger. 
(This has nothing to do with the corresponding inequality for $y_i$.)

We emphasize that this warning applies to 
nearly all of the following inequalities:
\begin{align}
 y_1y_2+y_1y_3+y_1y_4+y_2y_3 &\geq a_1a_2+a_1a_3+a_1a_4+a_2a_3\notag\\
 \vdots\\
 e_2^{(n)}(y)=y_1y_2+\cdots y_{n-1}y_n &\geq a_1a_2+\cdots a_{n-1}a_n=e_2^{(n)}(a)\,.\notag
\end{align}
Also for the products of three factors the estimate holds for the biggest, the sum of the two biggest, the sum of the biggest three, the biggest four...:
\begin{align*}
 y_1y_2y_3 &\geq a_1a_2a_3\notag\\
 y_1y_2y_3+y_1y_2y_4&\geq a_1a_2a_3+a_1a_2a_4\notag\\
 y_1y_2y_3+y_1y_2y_4+y_1y_3y_4&\geq a_1a_2a_3+a_1a_2a_4+a_1a_3a_4\\
 \vdots\notag\\
 e_3^{(n)}(y)&\geq e_3^{(n)}(a)\,,\notag
\end{align*}
and so on, until finally
\begin{align}
 y_1\,y_2\,\cdots\, y_n\,\;=\;\;e_n^{(n)}(y)&\;\;\geq\;\; e_n^{(n)}(a)\;\;=\;\;a_1\,a_2\,\cdots\, a_n\,.
\end{align}
{\bf Remark 1:} 
We arrived at the log-majorization in this paper by using the first condition each: the inequalities for $y_1$, $y_1y_2$, $y_1y_2y_3$ and so on.\\
The proof of the sum of squared logarithms inequality uses the last condition: $e_1(y)$, $e_2(y)$, $e_3(y)$ and so on.
\\[0.5cm]
{\bf Remark 2:} The sum of squared logarithms inequality is independent of Cohen's inequality \eqref{eq:Cohen}.\\[0.5cm]
{\bf Remark 3:} Cohen's theorem \eqref{eq:Cohen} can also be applied to the inverse matrices (as in \cite{NNF}).
The only additional inequality we gain is the other estimate for $e_n$, that is the equality of determinants. (Which we already know by different considerations.)
All the other ``new'' inequalities can be obtained by dividing the known ones by $y_1y_2\cdots y_n=a_1a_2\cdots a_n$.

\end{document}